\documentclass{article}

\usepackage{amsmath}

     \usepackage[final]{neurips_2024_ml4ps}

\usepackage[utf8]{inputenc} 
\usepackage[T1]{fontenc}    
\usepackage{hyperref}       
\usepackage{url}            
\usepackage{booktabs}       
\usepackage{amsfonts}       
\usepackage{nicefrac}       
\usepackage{microtype}      
\usepackage{xcolor}         
\usepackage{graphicx}
\usepackage{caption}
\usepackage{subcaption}
\newcommand{\RR}{I\!\!R} 

\newcommand{\bx}{\ensuremath{\boldsymbol{x}}}
\newcommand{\balpha}{\ensuremath{\boldsymbol{\alpha}}}
\newcommand{\bbeta}{\ensuremath{\boldsymbol{\beta}}}
\newcommand{\bxi}{\ensuremath{\boldsymbol{\xi}}}
\newcommand{\bmu}{\ensuremath{\boldsymbol{\mu}}}
\newcommand{\btheta}{\ensuremath{\boldsymbol{\theta}}}
\newcommand{\bpsi}{\ensuremath{\boldsymbol{\psi}}}
\newcommand{\calL}{\ensuremath{\mathcal{L}}}
\newcommand{\calC}{\ensuremath{\mathcal{C}}}
\newcommand{\calB}{\ensuremath{\mathcal{B}}}
\newcommand{\calD}{\ensuremath{\mathcal{D}}}
\newcommand{\calX}{\ensuremath{\mathcal{X}}}
\newcommand{\calT}{\ensuremath{\mathcal{T}}}

\usepackage{pifont}
\usepackage{enumitem}

\title{Physics-informed reduced order model with conditional neural fields}

\author{%
  Minji Kim \\
  Department of Statistics and Operations Research\\
  University of North Carolina at Chapel Hill\\
  \texttt{mkim5@unc.edu} \\
  \And
  Tianshu Wen \\
  Aerospace and Mechanical Engineering \\
  University of Notre Dame \\
  \texttt{twen2@nd.edu} \\
  \AND
  Kookjin Lee \\
  Computing and Augmented Intelligence \\
  Arizona State University \\
  \texttt{kookjin.lee@asu.edu} \\
  \And
  Youngsoo Choi \\
  Center for Applied Scientific Computing \\
  Lawrence Livermore National Laboratory \\
  \texttt{choi15@llnl.gov} \\
}

\begin{document}

\maketitle

\begin{abstract}
This study presents the conditional neural fields for reduced-order modeling (CNF-ROM) framework to approximate solutions of parametrized partial differential equations (PDEs). 
The approach combines a parametric neural ODE (PNODE) for modeling latent dynamics over time with a decoder that reconstructs PDE solutions from the corresponding latent states.
We introduce a physics-informed learning objective for CNF-ROM, which includes two key components. First, the framework uses coordinate-based neural networks to calculate and minimize PDE residuals by computing spatial derivatives via automatic differentiation and applying the chain rule for time derivatives. Second, exact initial and boundary conditions (IC/BC) are imposed using approximate distance functions (ADFs) [Sukumar and Srivastava, \textit{CMAME}, 2022]. However, ADFs introduce a trade-off as their second- or higher-order derivatives become unstable at the joining points of boundaries. To address this, we introduce an auxiliary network inspired by [Gladstone et al., \textit{NeurIPS ML4PS} workshop, 2022].
Our method is validated through parameter extrapolation and interpolation, temporal extrapolation, and comparisons with analytical solutions.
\end{abstract}

\section{Introduction}
Numerical simulations for solving nonlinear partial differential equations (PDEs) have advanced scientific understanding by enabling accurate modeling of complex physical phenomena. 
However, the computational demands of high-fidelity simulations have led to the development of various surrogate and reduced-order modeling (ROM) techniques. Some of these techniques simplify the underlying physics, while another approach is to accelerate computations by leveraging data, as seen in linear subspace projection-based ROMs \cite{CARLBERG2013623, kim2021efficient, copeland2022reduced, choi2019space, choi2020sns} and nonlinear manifold ROMs \cite{fulton2019latent,lee2020model,kim2022fast}.

In recent years, neural networks have been employed increasingly to approximate PDE solutions, with physics-informed neural networks (PINNs) \cite{RAISSI2019686,cai2021physics,cho2024parameterized} emerging as an approach that learns physics by incorporating the governing PDE directly into the loss function. A more recent direction has emerged in coordinate-based neural networks, originally developed to learn implicit representations of complex signals, commonly referred to as implicit neural representations (INRs) \cite{sitzmann2020implicit,tancik2020fourier}. 
These networks are particularly attractive for PDE-related tasks because they can learn continuous functions over domains. Building on this capability, 
a series of works has been proposed that use data-driven approaches to learn PDE solutions using an INR-based decoder \cite{yin2023continuous,chen2023crom,wan2023evolve}. This approach leverages the idea of modeling a vector field with a neural function that is conditioned on a latent state, a framework we refer to as conditional neural fields ROMs (CNF-ROMs). 

Building on these developments, our work extends the CNF-ROM framework to support both data-driven and physics-informed learning. 
To the best of our knowledge, this is the first application of a space- and time-separated CNF for training PINNs.
Our contributions are summarized as follows:
\begin{itemize}[leftmargin=*]
    
    \item Establishing a physics-informed learning objective for CNF-ROM framework,
    \item 
    Enabling the framework to handle parametrized PDEs using PNODEs \cite{lee2021parameterized,wen2023reduced}, 
    
    \item Addressing challenges in enforcing exact initial and boundary conditions (IC and BC) by introducing an auxiliary network to approximate first-order derivatives, following \cite{SUKUMAR2022114333,gladstone2023fopinns},
    
    \item Introducing training objectives with simultaneous optimization of decoder and PNODE parameters,
    \item Verifying the performance in parameter interpolation/extrapolation (with PINN fine-tuning for unseen parameters), temporal extrapolation, and comparison with analytical solutions. 
\end{itemize} 

\section{Models: Parameterized CNF-ROMs}
\label{s:setup}

We consider PDEs parametrized by $\bmu\in \calD$, with $\calD \in {\RR}^{N_{\bmu}}$, taking the form
\begin{align}
\label{eq:pde-govern}
    \partial_t u = \calL(u; \bmu),\quad u(\bx, 0, \bmu) = u_0(\bx, \bmu), \quad \calB(u;\bmu) = 0, \quad t\in [0,T],\quad \bx \in\Omega \in \RR^d,
\end{align}
where $\calL$ is a differential operator and $\calB$ is a boundary operator.
$u : \Omega \times (0,T]\times \calD \rightarrow {\RR}$ is the solution of the PDE that represents a physical quantity such as velocity or pressure. $u(\bx, t, \bmu)$ can be approximated by a neural network $u_{\btheta}(\bx, t, \bmu)$, where $\btheta\in \RR ^{d_{\btheta}}$ is a vector of tunable parameters.

\paragraph{Conditional neural fields}
\begin{figure}
    \centering
    \includegraphics[width=0.9\linewidth]{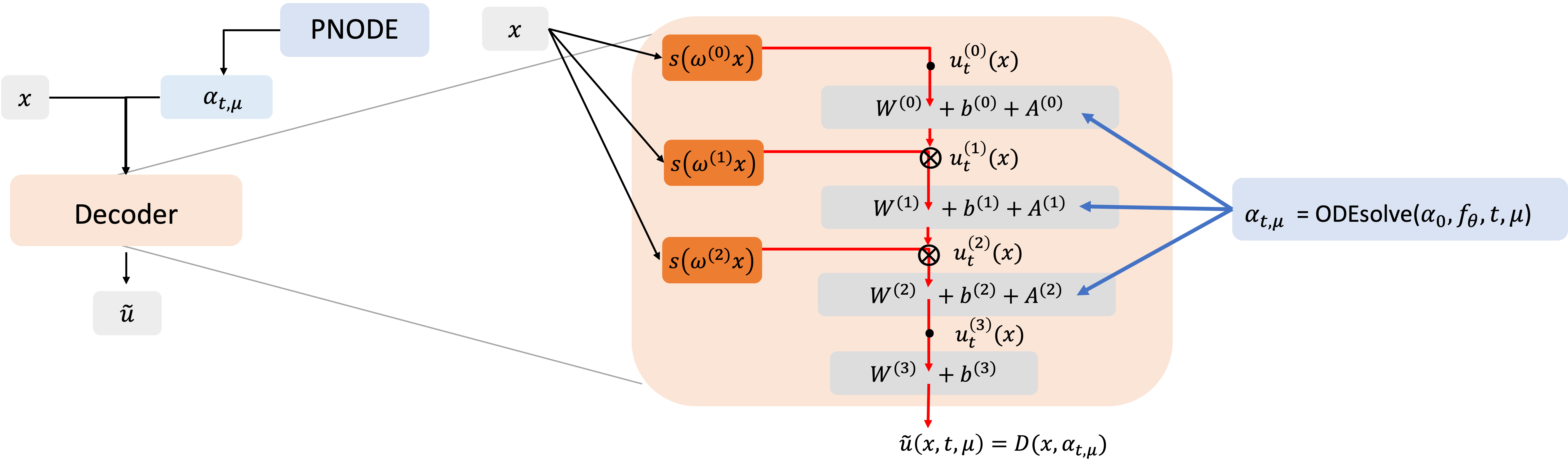}
    \caption{Diagram of the CNF-ROM structure, combining PNODE with the architecture of \cite{yin2023continuous}.}
    \label{fig:cnf-rom}
\end{figure}
When coordinate-based neural networks specifically operate in spatial and temporal domains, we refer to them as neural fields (NFs) \cite{NF_Xie2022}, borrowing the term ``fields" from physics.
CNFs \cite{pmlr-v235-hagnberger24a} extend NFs by introducing a conditioning latent factor $\balpha$.
A common approach uses a hypernetwork $h_{\bpsi}\!:\! \RR ^{d_{\balpha}}\!\!\rightarrow\!\! \RR ^{d_{\btheta}}$ to generate 
(a part of) high-dimensional parameters $\btheta = h_{\bpsi}(\balpha)\in \RR ^{d_{\btheta}}$ from a low-dimensional latent state $\balpha\in \RR ^{d_{\balpha}}$. 
This use of a hypernetwork enables reduced-order modeling by casting the problem as one of modeling latent states. We apply this framework to approximate solutions of spatio-temporal governing PDEs, referring to it as the conditional neural field for reduced-order modeling (CNF-ROM). This structure is illustrated in the diagram shown in Figure~\ref{fig:cnf-rom}.
The approach combines a parametric neural ODE (PNODE) \cite{lee2021parameterized}, which models latent dynamics, with a decoder that reconstructs PDE solutions. 
The decoder maps the latent state and coordinates to the solution, while the PNODE captures distinct latent trajectories for each $\bmu$.

\paragraph{Decoder} Building on the DINo architecture \cite{yin2023continuous}, we implement a decoder by employing a spatial coordinate-based network based on FourierNet \cite{fathony2021multiplicative}. To separate the spatial and temporal domains, DINo introduces the hypernetwork whose parameters are determined conditionally on the time-dependent latent state, which in our framework is extended to also depend on PDE parameters. 
The expanded section on the right side of Figure~\ref{fig:cnf-rom} illustrates the detailed DINo architecture integrated with PNODE.
Following this approach, we define our decoder as
\begin{align}
\label{eq:decoder}
\widetilde u(\bx, t, \bmu) = D_{{\bpsi}}(\bx, \balpha_{t, \bmu}),
\end{align}
where sinusoidal filters $s(\omega^{(l)} x) = [\sin (\omega^{(l)} x), \cos(\omega^{(l)} x)]^\top$ are used as Fourier basis.
At each time step $t$ and parameter $\bmu$, a latent state $\balpha_{t,\bmu}$ determines the high-dimensional parameters of the decoder, thereby shifting the learning of the decoder to the training of the hypernetwork parameters $\bpsi$. 


\paragraph{PNODE}
Recent advancements in PINNs emphasize the importance of learning parameterized solutions to generalize across multiple PDE parameter values ($\bmu$) \cite{cho2024parameterized}. 
These extensions aim to capture diverse solution behaviors as a function of parameters, avoid retraining for every new parameter, and identify patterns across the parameter space, enabling generalization to unseen parameter values.
Parameterized neural ODEs (PNODEs) \cite{lee2021parameterized} address these challenges by extending neural ODEs
(NODEs) \cite{chen2018neural} to incorporate $\bmu$ as inputs, allowing them to model the evolution of latent states for multiple trajectories. In our model, PNODEs learn the function $f_\theta$:
\begin{align}
\label{eq:pnode}
\frac{d \balpha_{t,\bmu}}{dt} = f_{\btheta}(\balpha_{t,\bmu},t, \bmu),
\end{align}
where $f_{\btheta}: \RR^{d_{\balpha}} \times [0,T]\times \calD \rightarrow \RR^{d_{\balpha}} $ represents the velocity of the latent state based on the current state, time, and PDE parameters.
By integrating PNODEs with the decoder in Eq.~\eqref{eq:decoder}, we model the evolution of parametrized PDE solutions over time in a reduced-order latent space. 

\section{Methods: Physics-informed training with exact satisfaction of IC/BC}
\label{s:method}

\paragraph{Physics-informed learning with parametrized CNF-ROM}
With the CNF-ROM framework introduced in Section~\ref{s:setup}, we propose a physics-informed learning objective for PINNs.
A physics-informed learning objective typically includes loss terms for the governing PDE residual, as well as for IC and BC. These terms ensure that the solution satisfies the PDE and adheres to the given IC and BC \cite{RAISSI2019686}.
To compute the necessary spatial and temporal derivatives for the PDE residual loss,
the CNF-ROM leverages coordinate-based neural networks to directly compute spatial derivatives via automatic differentiation. Temporal dynamics, on the other hand, are modeled through a latent state $\balpha_{t,\bmu}$ that evolves over time for each $\bmu$. Using the chain rule, temporal derivatives are expressed as: 
\begin{align}
\label{eq:chain-rule}
\partial_t D_{\bpsi}(\bx, \balpha_{t,\bmu}) = \partial_{\balpha} D_{\bpsi}(\bx, \balpha_{t,\bmu}) \cdot f_{\btheta}(\balpha_{t,\bmu}, t, \bmu),
\end{align}
where $\partial_{\balpha} D_{\bpsi}(\bx, \balpha_{t,\bmu}) = \left.\partial_{\balpha} D_{\bpsi}(\bx, \balpha)\right|_{\balpha_{t,\bmu}}$ can also be computed using automatic differentiation.
This combination of coordinate-based neural networks and latent-state dynamics enables efficient and accurate computation of spatio-temporal derivatives.

Instead of incorporating IC and BC as loss terms, we address strategies for imposing exact IC and BC and discuss the trade-offs in this approach in the following paragraphs.

\paragraph{Exact imposition of IC and BC}
Exact imposition of IC and BC reduces the number of loss constraints while ensuring the uniqueness of the solution, which improves training convergence and predictive accuracy. 
The R-function-based approximate distance functions (ADFs) $\phi$ proposed by \cite{SUKUMAR2022114333} were developed to enforce boundary conditions while satisfying other desirable properties such as differentiability. These functions can be constructed to satisfy Dirichlet, Neumann, and Robin conditions a priori, even on complex and irregular geometries; see \cite{SUKUMAR2022114333} for more details.
Building on this foundation, we extend ADFs to the temporal domain to address initial conditions as well by ensuring \(\phi(\bx, t) = 0\) when \(\bx \in \partial \Omega\) or \(t = 0\).
For Dirichlet boundary conditions, where $\widehat u = g$ is imposed on boundaries, we use the following construction:
\begin{align}
\label{eq:exact-u}
\widehat{u}(\bx,t, \bmu;\bpsi, \btheta) = g(\bx,t;\bmu) + \phi(\bx, t) D_{\bpsi}(\bx, \balpha_{t,\bmu}(\btheta)).
\end{align}
Here, $\phi$ is an ADF and $D_{\bpsi}(\bx, \balpha_{t,\bmu})$ denotes the decoder output in Eq.~\eqref{eq:decoder}.

When using the R-function-based ADFs, the following challenges arise when the input dimension is greater than two. In our case, this holds as the ADF includes both time and space coordinates.
To be specific, (i) ADF $\phi$ is not well-defined on the boundary, and (ii) its second (or higher) derivatives explode near boundary junctions, as noted in \cite{SUKUMAR2022114333, gladstone2023fopinns, BERRONE2023e18820}. To mitigate these challenges, (i) we exclude the boundary set from the PINN loss definition in Eq.~\eqref{eq:loss-pde1} and (ii) adopt a more delicate approach proposed in \cite{gladstone2023fopinns}.
This approach addresses the trade-off in imposing IC and BC on $\widehat u$ by approximating its first-order derivatives. For example, to compute $\partial_{xx} \widehat{u}$, instead of directly computing the second-order derivatives, we first train an auxiliary network to approximate the first-order derivatives and then compute the derivative of this network’s output. To achieve this, we train another CNF-ROM $\widetilde v(\bx, t, \bmu; \bxi, \btheta) = D_{\bxi}(\bx,\bbeta_{t,\bmu}(\btheta))$ and add an additional loss term matching the derivatives,
\begin{align}
\label{eq:loss-deriv}
L_{\text{deriv}}(\bxi, \btheta) = \frac{1}{ N_o}\sum_{\bmu \in \calD}\sum_{t \in \calT\backslash \{0\}^c}\sum_{\bx \in \calX\backslash \partial \Omega^c} \| \partial_x \widehat u (\bx,t,\bmu)- \widetilde v(\bx,t,\bmu;\bxi, \btheta) \|^2,
\end{align}
so that the second or higher-order derivatives of $\widehat u$ are approximated using the derivatives of $\widetilde v$.
Note that $\widetilde v$ is constructed without ADF constraints, and thus higher-order derivatives are obtainable.

In the following paragraph, we describe two learning modes: data-driven learning and physics-informed learning. In both modes, the solution is constructed as in Eq.~\eqref{eq:exact-u} to impose IC and BC.

\paragraph{Training objectives}

One key advantage of imposing exact IC and BC is that it eliminates the need for an encoder. Since the output satisfies the initial condition for any arbitrary \(\balpha_0\), we initialize \(\balpha_0 = \boldsymbol{0}\in \RR^{d_{\balpha}}\). Once the initial latent state is obtained, we use ODE solvers \cite{chen2018neural} to evolve it according to Eq.~\eqref{eq:pnode}. 
The use of PNODEs allows the initial trajectory to adapt based on different values of \(\bmu\), enabling the model to generate distinct outputs.
Note that for $\widetilde v$ in Eq.~\eqref{eq:loss-deriv}, $\bbeta_{0,\bmu}$ needs to be learned, and we use an auto-decoding approach as suggested in \cite{Park_2019_CVPR}.

We propose the following training objectives for CNF-ROM, employing an efficient simultaneous optimization that jointly updates the decoder and PNODE parameters.
\begin{itemize}[leftmargin=*]
\item {\it Data loss:} When data is available, we minimize the following data loss:
\begin{align}
\label{eq:loss-data}
L_{\text{data}}(\bpsi, \btheta) =\frac{1}{N} \sum_{(\bx, t, \bmu) \in \calC}\| u(\bx, t, \bmu) - \widehat u(\bx,t,\bmu; \bpsi, \btheta) \|^2,
\end{align}
where $\widehat u$ is defined in Eq.~\eqref{eq:exact-u}, $\calC = \calX\times \calT\times \calD$ is the set of collocation points, and $N=|\calC|$.
\item {\it PINN loss:} For physics-informed training, we minimize the following PDE residual loss:
\begin{align}
\label{eq:loss-pde1}
 L_{\text{PDE}}(\bpsi, \btheta) = \frac{1}{N_o}\sum_{(\bx, t, \bmu) \in \calC_0} \| \partial_t \widehat u (\bx,t,\bmu; \bpsi, \btheta)- \calL(\widehat u(\bx,t,\bmu; \bpsi, \btheta); v(\bx, t, \bmu)) \|^2,
\end{align}
where $\calC_0 = \calX\backslash \partial \Omega^c \times \calT\backslash \{0\}^c\times \calD , ~N_o = |\calC_0|,$ and $\calL$ is from Eq.~\eqref{eq:pde-govern}. When the governing equation involves second or higher-order derivatives, these terms in Eq.~\eqref{eq:loss-pde1} are replaced with those computed using the derivatives of $\widetilde v$. Note that IC and BC are imposed a priori.
\end{itemize}

We are now equipped with the training losses defined in Eq.~\ref{eq:loss-deriv}, \eqref{eq:loss-data} and \eqref{eq:loss-pde1}. In Section~\ref{s:result}, we present training scenarios and numerical results for the 1D viscous Burgers equation.

\section{Numerical results and discussion}
\label{s:result}
We consider the 1D viscous Burgers equation on $\Omega \times [0,T] = [0,2] \times [0,1]$ with specific IC and BC, where a known solution exists. The governing equation, parametrized by the Reynolds number $\mu$, is: 
\begin{align}
    \partial_t u + u \cdot \partial_x u - (1/\mu) \partial_{xx} u = 0, \quad u(0,t) = u(2,t) = 0 \text{ on } t \in [0, 1].
\end{align}
The exact solution used to validate the performance of our models is given by
\begin{align}
\label{eq:cecchi}
    u_{ex}(x,t;\mu) = \frac{2\pi}{\mu} \left\{\frac{\frac{1}{4}e^{-\pi^2 t/\mu}\sin(\pi x) + e^{-4\pi^2 t/\mu}\sin(2\pi x)}{1 + \frac{1}{4}e^{-\pi^2 t/\mu}\cos(\pi x) + \frac{1}{2}e^{-4\pi^2 t/\mu}\cos(2\pi x)} \right\},
\end{align}
with the initial condition determined by \eqref{eq:cecchi} at $t=0$. 
The sets of training and test parameters $\bmu_{\text{train}} = \{20, 30, 40, 50, 60, 70, 80, 100\}$ and $\bmu_{\text{test}} = \{15, 25, 45, 90, 110\}$ and the latent state dimension $d_{\balpha}=10$ are used. The numerical solution of the full-order model was obtained using a uniform spatial and temporal discretization (\(n_x=64\), \(n_t=100\)) and a backward Euler time integrator.
We consider the following training scenarios:
\begin{itemize}
\item[(a)] Training with data loss: $L_{\text{data}}(\bpsi, \btheta) + L_{\text{deriv}}(\bxi, \btheta)$, $T=1$, across $\bmu_{\text{train}}$.
\item[(b)] Fine-tuning with PINN loss: $L_{\text{PDE}}(\btheta) + L_{\text{deriv}}(\btheta)$, $T=1$, $\bmu \in \bmu_{\text{train}} \cup \bmu_{\text{test}}$.
\end{itemize}
The goal of these scenarios is to highlight the advantage of using PINN as a fine-tuning objective, leveraging its ability to learn without data.
For scenario (a), both the decoder and the auxiliary network are optimized during training. For scenario (b), starting from the pre-trained model from (a), we freeze the decoder and use the PDE loss to update only the PNODE (latent dynamics) parameters for the target parameters. This PINN fine-tuning stage models only the low-dimensional latent state, aligning with the ROM perspective.

\begin{figure}[t]
\centering
\begin{subfigure}[b]{0.58\textwidth}
\centering\includegraphics[width=\textwidth,height=2.8cm]{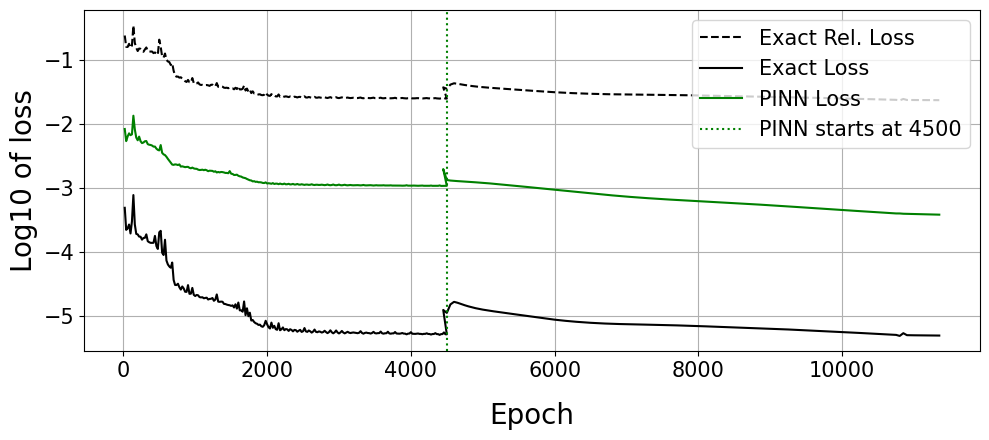}
\vspace{0.2cm} 
\includegraphics[width=\textwidth,height=5cm]{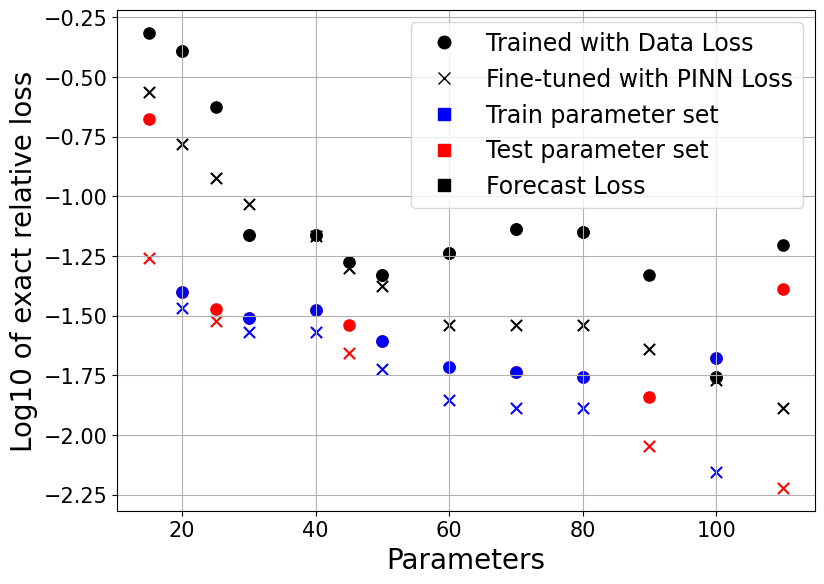}
\end{subfigure}
\hfill
\begin{subfigure}[b]{0.38\textwidth}
\centering
\includegraphics[width=\textwidth,height=3.9cm]{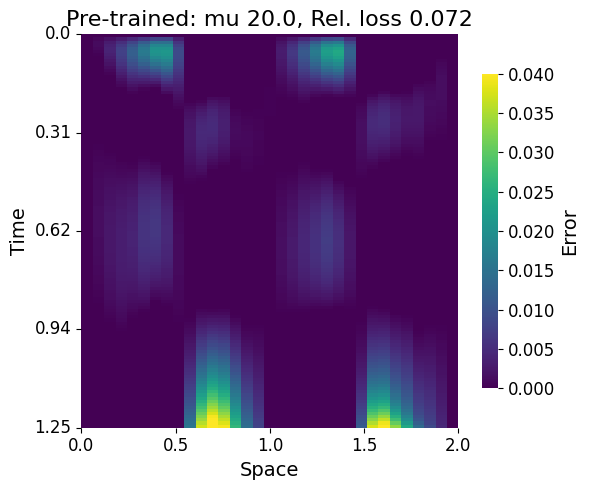}
\vspace{0.2cm}
\includegraphics[width=\textwidth,height=3.9cm]{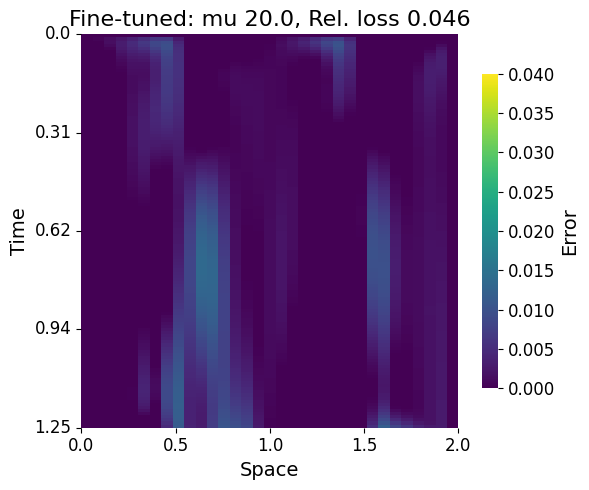}
\end{subfigure}
\caption{Performance verification through: (top left) loss trajectories, (bottom left) evaluation on training and test set parameters, and temporal extrapolation at $t=1.25$; (top right) heatmap of loss for model (a) at $\bmu = 20$; (bottom right) heatmap of loss for model (b) at $\bmu = 20$.}
\label{f:validation}
\end{figure}

Figure~\ref{f:validation} presents the performance evaluation from multiple perspectives. 
The top left panel of presents the loss trajectory, displaying the loss with respect to exact solutions in Eq.~\eqref{eq:cecchi} (``Exact Loss", solid black), its relative error (``Exact Rel. Loss", dashed black) and the PDE loss for PINN (``PINN Loss", solid green). The model was pre-trained under scenario (a) up to epoch 4500, followed by fine-tuning under scenario (b) for the fixed parameter $\bmu = 20$. 
The purpose of this figure is to validate the PDE loss calculation by showing that both exact and PINN losses decrease in parallel, regardless of whether the training was based on (a) data loss or (b) PINN loss, confirming that the PDE residual loss aligns with the exact solution.
The right panels of Figure~\ref{f:validation} present heatmaps of the loss for \(\mu = 20\). The top right heatmap corresponds to the model trained under scenario (a), while the bottom right heatmap shows the result of the fine-tuned model (b). We first observe that additional fine-tuning with PINN loss led to additional error reduction. Notably, the region \(t > 1\) represents an extrapolation beyond the trained domain (i.e. forecasting). The error plots reveal that forecast errors are more pronounced in the pre-trained model, whereas the fine-tuned PINN loss model exhibits robust result in this region.

The bottom left panel of Figure~\ref{f:validation} offers a summary of the performance of the model trained under the pre-trained scenario (a, \ding{108}) and the fine-tuned scenario (b, \ding{53}) for each parameter \(\bmu \in \bmu_{\text{train}} \cup \bmu_{\text{test}}\). This panel compares fine-tuning results for the training set (blue), inter-/extrapolation results in the parameter space (red), and temporal extrapolation (black).
For temporal extrapolation, all models (a) and (b) were trained up to $T=1$ and used to evaluate up to $t=1.25$. Specifically, fine-tuning with PINN loss demonstrated clear improvements across all regions except temporal extrapolation performance at $\mu=30$, with a particularly significant error reduction in parameter extrapolation areas where $\mu=15$ or $110$.
In general, we observed that the fine-tuning with PINN loss improves model performance by learning without additional data, enabling further learning for unseen parameters, and enhancing robustness in forecast regions.
Finally, we conclude with the remark that the CNF-ROM framework can be extended to higher-dimensional problems, which we leave their exploration as a future direction.



\begin{ack}
\textbf{K. Lee} was supported by NSF under grant IIS \#2338909.
\textbf{Y.\ Choi} was supported by the U.S. Department of Energy, Office of Science, Office of Advanced Scientific Computing Research, as part of the CHaRMNET Mathematical Multifaceted Integrated Capability Center (MMICC) program, under Award Number DE-SC0023164.
Lawrence Livermore National Laboratory is operated by Lawrence Livermore National Security, LLC, for the U.S. Department of Energy, National Nuclear Security Administration under Contract DE-AC52-07NA27344. IM release number: LLNL-CONF-869137.
\end{ack}

{
\small
\bibliography{main}
\bibliographystyle{unsrt}
}






\end{document}